\documentclass[11pt,a4paper]{amsart}
\usepackage[utf8]{inputenc}
\usepackage{amsmath}
\usepackage{amssymb}
\usepackage{amsthm}
\usepackage{stmaryrd}
\usepackage{graphicx}
\usepackage{tikz}
\usepackage{tikz-cd}
\usepackage{mathtools}
\usepackage{comment}
\usepackage{enumitem}
\usepackage{import}

\usepackage[
bookmarks         = true,
bookmarksnumbered = true,
pdfstartview      = FitH,
pdfpagelayout     = OneColumn, 
pdfpagemode=UseOutlines,
colorlinks        = true,
linkcolor         = blue,
urlcolor          = magenta,
citecolor         =red,
backref           = page
]{hyperref}

\usepackage[draft,inline,nomargin,index]{fixme}

\fxsetup{theme=color,mode=multiuser}
\FXRegisterAuthor{ad}{aad}{\color{purple}AD}
\FXRegisterAuthor{jl}{aajl}{\color{blue}JL}

\usepackage{fourier-orns}
\usepackage{natbib}
\usepackage{ mathrsfs }
\DeclareMathOperator{\Aut}{Aut}

\DeclareMathOperator{\pr}{pr}

\DeclareMathOperator{\supp}{supp}

\renewcommand{\d}{\mathrm d}

\renewcommand{\le}{\leqslant}
\renewcommand{\ge}{\geqslant}
\renewcommand{\epsilon}{\varepsilon}
\renewcommand{\phi}{\varphi}
\renewcommand{\rm}{\mathrm}

\renewcommand{\L}{\mathrm L}

\newcommand{\m}{\mathbf}
\newcommand{\tendv}{\xrightarrow[n\rightarrow +\infty]{}}

\newcommand{\mc}{\mathcal}

\newcommand{\inv}{^{-1}}
\newcommand{\Prob}{\rm{Prob}}

\theoremstyle{definition}
\newtheorem{dfn}{Definition}[section]

\theoremstyle{plain}
\newtheorem{theo}[dfn]{Theorem} 
\newtheorem{pro}[dfn]{Proposition}
\newtheorem{lem}[dfn]{Lemma}
\newtheorem{cor}[dfn]{Corollary}
\newtheorem{prodfn}[dfn]{Proposition-Definition}

\newtheorem{theointro}{Theorem}

\newtheorem{corintro}[theointro]{Corollary}

\theoremstyle{remark}
\newtheorem{rem}[dfn]{Remark}
\newtheorem{ex}[dfn]{Example}

\title[Boundaries of lattices]{Poisson boundaries of lattices in locally compact groups}
\author{Antoine Derimay}
\date{}
\begin{document}
\begin{abstract}
    For a random walk on a general locally compact group we provide a discretization procedure transforming the walk to one on a lattice of the group, preserving the Poisson boundary and some finiteness of moments.
\end{abstract}
\maketitle

\tableofcontents 
\section{Introduction}
The Poisson-Furstenberg boundary is a powerful measure-theoretic object associated to a general Markov chain, which keeps track of the asymptotic behaviour of the chain as well as bounded harmonic functions for that chain. If the Markov chain happens to be a right random walk on a group, say $G$ with law $\mu$, the Poisson boundary inherits a measure-class preserving $G$-action, which is both amenable in the sense of Zimmer and ergodic. Moreover, it is $\mu$-stationary, and it is in a sense the maximal such object, all $\mu$-stationary $G$-spaces being (measure-preserving extensions of) a factor of the Poisson boundary. Finally, if the step measure $\mu$ is symmetric, the Poisson boundary is the prime example of a strong boundary in the sense of Bader-Furman \cite{BF14}, opening up the door for a number of applications to rigidity theory if its combinatorics is well understood \cite{BF14,BCFS24,GHL22,Der25}.\medskip

Hence the Poisson boundary is a versatile object, whose understanding can be used in a number of contexts and applications. Amenability of the $G$-action in particular means that for any compact metrizable space $X$ on which $G$ acts, there is a $G$-equivariant map from the Poisson boundary of $G$ to the set of probability measures on $G$. What Furstenberg noticed in \cite{Fur73} is that this property passes to lattices in $G$: if $\Gamma$ is a lattice in $G$ and $\mu$ is a ``reasonable" measure on $G$, then for any compact metrizable $\Gamma$-space $X$ there is a $\Gamma$-map from the Poisson boundary of $G$ to the space of probability measures on $X$. This fact then turns out to be a crucial ingredient of the celebrated Margulis superrigidity and normal subgroup theorems (see \cite{Mar91}, the discretization procedure is Proposition VI.4.1).

The proof of Furstenberg's result hinges on a ``discretization procedure", the idea of which has then been improved several times in the context of Brownian motion on Riemannian manifolds, first by Lyons and Sullivan in \cite{LS84}, see also \cite{Kai92,BL96}, where it is proven that for a big enough discrete subset of the manifold, there is a Markov chain on this subset for which the restrictions of bounded harmonic functions on the manifold are harmonic, and the datum of this restriction is enough to recover the whole function. In other words, the Poisson boundaries of the manifold and the Markov chain are naturally isomorphic.

This works very well for lattices in Lie groups, where this construction has enabled a number of rigidity results, including the recent paper \cite{BCFS24}, in which it is proven that lattices in semisimple Lie groups may not act on hyperbolic spaces that are not the symmetric spaces of some rank-one factor (see also \cite{DH20,ALSW26} for instance).

It works much less well in situations where Brownian motion cannot be defined, for instance in totally disconnected groups. The solution then becomes an adaptation of the Lyons-Sullivan construction to processes with discrete time increments but non-discrete state space, much as in \cite{KM96} where random walks on Teichmüller space are considered.

This adds a number of technical complications however, the biggest of them is that the \emph{Harnack inequality} is not necessarily true for such systems, and is an assumption we have to add, and which we call here \emph{condition (H)} (see Definition \ref{dfn: condition (H)}). Under this condition we are able to give a new discretization procedure:

\begin{theointro}\label{theointro: LS pour des marches sur des groupes lcsc}
    Let $G$ be a locally compact second countable group, and let $\Gamma$ be a lattice in $G$. Let $\mu\in \Prob(G)$ be generating, and assume moreover that $\mu$ satisfies condition (H) (see Definition \ref{dfn: condition (H)}). Then there exists a probability measure $\nu$ on $\Gamma$ such that 
    $$B_{PF}(\Gamma,\nu)\simeq B_{PF}(G,\mu).$$
    Moreover, $\nu$ has finite first moment as soon as $\mu$ does, and if additionally $\Gamma$ is uniform, $\nu$ has finite $k$-th moment as soon as $\mu$ does, and has a finite exponential moment as soon as $\mu$ does.
\end{theointro}
Additionally the law of the random walk on $\Gamma$ is well approximated by a sequence of stopping times defined on the random walk on $G$ (much like in the original Lyons-Sullivan procedure); this is made precise in Proposition \ref{pro: approximation by N_k}\medskip

Condition (H) is a strong condition, not verified for many measures, but measures satisfying (H) may be constructed from nicely behaved measures without changing the Poisson boundary (see Proposition \ref{pro: pullback de RW de graphe a (H)} and Corollary \ref{cor: supp cpt-> mu chapeau a (H)} for example). The following is Theorem \ref{theointro: LS pour des marches sur des groupes lcsc}, along with Proposition \ref{pro: mu chapeau conserve la frontière de Poisson} and Corollary \ref{cor: supp cpt-> mu chapeau a (H)}:

\begin{corintro}\label{corintro: discrétisation sans parler de (H)}
    Let $G$ be a locally compact second countable group with Haar measure $\lambda$, and let $\Gamma$ be a lattice in $G$. Let $\mu\in \Prob(G)$ be admissible and such that $\mu\ll\lambda$ with $\d\mu/\d\lambda$ continuous and of compact support. Then there exists a probability measure $\nu$ on $\Gamma$ such that 
    $$B_{PF}(\Gamma,\nu)\simeq B_{PF}(G,\mu).$$
    Moreover, $\nu$ has finite first moment as soon as $\mu$ does and if additionally $\Gamma$ is uniform, $\nu$ has finite $k$-th moment as soon as $\mu$ does, and has a finite exponential moment as soon as $\mu$ does.
\end{corintro}

One thing Theorem \ref{theointro: LS pour des marches sur des groupes lcsc} allows us to do is to get good descriptions of (possibly irreducible) lattices in products of groups:

\begin{corintro}\label{corintro: BPF de réseau dans produit}
    Let $G_i,1\le i\le n$ be lcsc groups, and $\mu_i$ be generating probability measures on $G_i$ with property (H). Let $\Gamma$ be a lattice in the product $G=\prod G_i$. Then there exists a generating probability measure $\nu$ on $\Gamma$ such that 
    $$B_{PF}(\Gamma,\nu)=B_{PF}(G_1,\mu_1)\times\cdots\times B_{PF}(G_n,\mu_n).$$
\end{corintro}

\begin{rem}
    This proves in particular that an $S$-arithmetic lattice $\Gamma=\underline G(\mc O_S)$ admits the product of the $\underline G(K_v)/\underline P(K_v)$ for $v\in S$ as a Poisson boundary (see \cite{Mor15,PR94} for more information on $S$-arithmetic lattices).\medskip

    This does \emph{not} however answer Remark 1.6 of \cite{BCFS24}, since so far there is no guarantee that the step measure on $\Gamma$ can be chosen to be symmetric. This will be done using another discretization procedure in a more restricted setting in an upcoming project.
\end{rem}

Finally, we prove a generalization of \cite[Theorem 3.10]{Der25} to unimodular locally compact groups acting on a graph with finite covolume, which is Theorem \ref{theo: RW on graph -> on group}. Besides being useful as is, it could be useful while having Corollary \ref{corintro: BPF de réseau dans produit} in mind in order to determine a Poisson boundary of a lattice in a product of groups, where some act on discrete spaces. \medskip

\noindent\textbf{Organization of the paper.}
Section \ref{sec: Random walks} is devoted to the presentation of probabilistic background notions that will be useful later on: Markov chains on non-discrete state spaces, recurrence properties and Poisson boundaries. Section \ref{sec: discretisation} is devoted to proving Theorem \ref{theointro: LS pour des marches sur des groupes lcsc} and the related results. Finally, in Section \ref{sec: kaima woess+recurrence} we establish Theorem \ref{theo: RW on graph -> on group} and linked results.\medskip

\noindent\textbf{AI statement.} no AI whatsoever has been used in this project.\medskip

\noindent\textbf{Acknowledgements.}
The author would like to thank Jean Lécureux for numerous comments and advice on previous drafts of this paper and Alano Ancona for an enjoyable discussion about discretisation procedures. Additionally, he thanks Ebrahim Samei for showing interest in this work, and pushing him to publish it.
\newpage

\section{Random walks}\label{sec: Random walks}

\subsection{Discrete-time Markov chains}
We present here what we need of the theory of Markov chains for a continuous state space, but with discrete increments in time. For this part, we follow the book \cite{Rev84}.\medskip

\subsubsection{Markov kernels}
The typical image to have in mind here is that of a particle on $\m R^d$, sampled on regular time intervals, say every second. Then according to where the particle is at time $n$, it samples a normal distribution (independently from whatever happened at previous times), and at time $n+1$ it has moved by that amount. 

\begin{dfn}
    Let $(X,\mc E)$ be a measurable space. A \emph{Markov kernel} on $(X,\mc E)$ is a map 
    $$P:X\times \mc E\to [0,1]$$
    such that for every $x\in X$, the map $A\mapsto P(x,A)$ is a probability measure on $X$ (which we will sometimes denote $P(x,\cdot)$), and for every $A\in \mc E$, the map $x\mapsto P(x,A)$ is measurable.
\end{dfn}

A Markov kernel $P$ should be thought of as a way of encoding all the probability distributions $P(x,\cdot)$ that, given that the particle is at $x$ at time $n$, it is with probability $P(x,A)$ in $A$ at time $n+1$. Denote by $Z_n$ the random variable which associates the position of the particle at time $n$. We will construct a precise model for the $Z_n$ later on.

\begin{ex}\label{ex: Markov kernels}\:
    \begin{enumerate}[label=(\roman*)]
        \item \label{ex: Markov kernels; discrete}If $X$ is countable and $\mc E$ is the discrete $\sigma$-algebra, it suffices to know $P$ on elements of the form $P(x,\{y\})$, which we may as well write $P(x,y)$. In that case, we recover the notion of a random walk on a discrete space.
        \item \label{ex: Markov kernels; on groups}If $G$ is a lcsc topological group and $\mu$ is a probability measure on $G$, we may define a kernel $P$ by the formula 
        $$P(x,A)=\mu(x\inv A).$$
        This is called the \emph{right} random walk on $G$ with law $\mu$ (the left random walk has kernel $P(x,A)=\mu(Ax\inv)$). One usually considers the \emph{right} random walk, because the natural objects associated to it (path space, Poisson boundary...) admit a \emph{left} $G$-action.
        \item \label{ex: Markov kernels; with group action}More generally, if $G$ is a lcsc topological group which acts continuously on the left on some topological space $M$ and $\mu$ is a probability measure on $G$, we also have a kernel $P$ on $M$ defined by 
        $$P(x,A)=\mu(\{g\in G,gx\in A\}).$$
    \end{enumerate}
\end{ex}

To a kernel $P$ one can also associate an action on functions to $\m R$. Given a measurable map $f:X\to \m R$, one can let 
$$Pf(x):=\int_XP(x,\d y)f(y)=\int_Xf(y)\d P(x,\cdot)(y),$$
under the assumption that the integral makes sense. It is the average of the values of $f$, weighted by the distribution of $Z_1$ under $Z_0=x$. A function $f$ such that $Pf=f$ is said to be $P$\emph{-harmonic} (or simply harmonic when $P$ is clear from context), and ($P$-)\emph{superharmonic} when $Pf\le f$. The constant maps are always harmonic.

Dually, $P$ also acts on positive measures by 
$$mP(A):=\int_X P(x,A)\d m(x)$$
(the convention is to write $P$ to the right of $m$ in this case). We will say $m$ is \emph{stationary} if $mP=m$. 

The operator $P$ is self-adjoint, in the sense that $mP(f)=m(Pf)$, where $\mu(f)=\int_Xf\d\mu$.\medskip

Given $P,Q$ two Markov kernels on $(X,\mc E)$, their \emph{composition} (or product) is the kernel $PQ$ defined by
$$PQ(x,A)=\int_X P(x,\d y)Q(y,A).$$
Composition can be verified to be associative. We will denote $PP$ by $P^2$, and more generally $P^n=PP^{n-1}=P^{n-1}P$. The measure $P^n(x,\cdot)$ is the distribution of the position of the particle at time $n$ given that it was at $x$ at time $0$.\medskip

\subsubsection{The canonical Markov chain}

We now construct explicitly the variables $Z_n$ which behave like the time $n$ map of the moving particle.

\begin{prodfn}\label{prodfn: canonical Markov chain}
    For any $x\in X$, there is a unique probability measure $\m P_x$ on $(X^\m N,\mc F):=(X,\mc E)^\m N$ such that for any $n\in \m N$ and any $A_0,A_1,\cdots,A_n\in \mc E$, 
    \begin{eqnarray*}
        &&\m P_x\left(A_0\times A_1\times\cdots\times A_n\times X^{\m N_{>n}}\right)\\
        &=&\m1_{A_0}(x)\int_{A_1}P(x,\d x_1)\int_{A_2}P(x_1,\d x_2)\cdots\int_{A_{n-1}}P(x_{n-2},\d x_{n-1})P(x_{n-1},A_n).
    \end{eqnarray*}
    For $n\in\m N$,define the random variable $Z_n:X^\m N\to X$ by $Z_n((x_0,x_1,\ldots))=x_n$. The sequence $(Z_n)_n$ is then called the \emph{canonical Markov chain} with transition probability $P$.
\end{prodfn}

\begin{proof}
    See \cite[Theorem 2.8]{Rev84}.
\end{proof}

More generally, given a measure $\nu$ on $X$, one defines a measure $\m P_\nu$ on $\Omega$ by 
$$\m P_\nu(A)=\int_X\m P_x(A)\d\nu(x).$$
It is a probability whenever $\nu$ is. In that case it is the distribution of the path of the particle, given that its starting distribution was $\nu$. \medskip

Our previously-defined notions now have useful probabilistic interpretations. For $n\in\m N$, let $\mc F_n$ be the $\sigma$-algebra generated by the $Z_m,m\le n$. It is the $\sigma$-algebra which ``remembers" the position of the particle at times up to $n$. Denote also by $\m E_x$ the expectation associated to $\m P_x$. 

The measure $P(x,\cdot)$ is (more or less by definition) the distribution of $Z_1$ with respect to $\m P_x$, and more generally $P^n(x,\cdot)$ is the distribution of $Z_n$ with respect to $\m P_x$. A map $f:X\to\m R$ is harmonic exactly when 
$$\m E_x[f(Z_n)]=f(x),$$
and equivalently when 
$$\m E_x[f(Z_{n+1})|\mc F_n]=f(Z_n),$$
in other words when $(f(Z_n))_n$ is a martingale with respect to the filtration $\mc F_n$. Similarly, $f$ is superharmonic exactly when $(f(Z_n))_n$ is a supermartingale.

Finally, the measure $mP$ is the distribution of $Z_1$ with respect to $\m P_m$ (and the distribution of $Z_n$ with respect to $\m P_m$ is $mP^n$).

\subsubsection{Stopping times and balayage}
A \emph{stopping time} $T$ with respect to the filtration $\mc F_n$ is a random variable $T:\Omega\to\m N$, such that for any $n\in\m N$, $\{T=n\}\in\mc F_n$. 

\begin{ex}
    The stopping times that we will encounter most of the time are what are called \emph{first hitting times}. If $A\in\mc E$, we define 
    $$T_A=\inf\{n\ge 0, Z_n\in A\}$$
    to be the first hitting time of $A$. It is easily checked to be a stopping time.
\end{ex}

If we assume further that $T_A$ is finite a.s. (compare with Definition \ref{dfn: Harris recurrent set}), let $Z^A:=Z_{T_A}$ be defined by $Z^A(\omega)=Z_{T_A(\omega)}(\omega)$. We may then define the \emph{balayage} transformation in the following way: for $x\in X$, let $\delta_x^A$ be the measure defined by
$$\delta_x^A(B)=\m P_x(Z^A\in B).$$
Since $T_A$ is finite a.s., $\delta_x^A$ is a probability measure, and it is supported on $A$.

Now for any positive measure $\nu$ on $X$, let 
$$\nu^A:=\int_X\delta_x^A\d\nu(x),$$
or equivalently, $\nu^A(B)=\m P_\nu(Z^A\in B)$. The measure $\nu^A$ is called the \emph{balayage measure} of $\nu$ with respect to $A$. It has same mass as $\nu$, and is supported on $A$. The properties of balayage measures we will be using later on are the following:

\begin{lem}\label{lem: properties of balayage}\:
    \begin{enumerate}
        \item The map $\nu\mapsto \nu^A$ is additive and homogeneous on the cone of positive measures on $X$.
        \item For any bounded harmonic function $h$, $\nu(h)=\nu^A(h).$
    \end{enumerate}
\end{lem}

\begin{proof}
    The first point follows from the definition. The second is a consequence of Doob's optional stopping time theorem: the sequence $h(Z_n)$ is a bounded martingale and $T_A$ is finite a.s., so
    $$\nu^A(h)=\m E_\nu(h(Z_{T_A}))=\m E_\nu(h(Z_0))=\nu(h).$$
\end{proof}

\subsubsection{The case of groups}
The special case of Example \ref{ex: Markov kernels} \ref{ex: Markov kernels; on groups} will be of a particular interest to us. We review here how the notions described above interact with notions that may be defined in the group-theoretic setting. Fix for the rest of this section a topological lcsc group $G$ and a measure $\mu$ on $G$. Let $P$ be the kernel associated with the right random walk on $G$ with law $\mu$.

A map $f$ will be ($\mu$-) harmonic exactly if 
$$f(g)=\int_Gf(gh)\d\mu(h).$$

Keeping the same notation, we also have 
$$mP(A)=m*\mu(A)=\int_G\int_G\m 1_A(gh)\d m(g)\d\mu(h)$$
for a measure $m$ on $G$. In particular, if there is a measure that is invariant under the $G$-action, it will automatically be stationary.

The interest of studying groups is the additional presence of a group action of $G$ on the path space $X^\m N=G^\m N$, being the diagonal left action. It preserves the walk, in the sense that 
$$P(gx,gA)=P(x,A).$$
Finally, from all this we see that we also have $g_*\m P_x=\m P_{gx}$, and more generally $g_*\m P_\nu=\m P_{g_*\nu}$.

\subsection{Irreducibility, recurrence}
This section is mostly taken from \cite{MT09} and \cite{Rev84}, and much more details can be found in either reference.

Recurrence means that the walk "comes back" infinitely often.  In the discrete setting, recurrence of the walk is only considered on \emph{irreducible} chains, that is, chains for which any point has a positive probability of visiting any other point. This is sensible, since it will be hard to reach a set with probability $1$ if you know that you reach it with probability $0$. 

If the measures $P(x,\cdot)$ are non-atomic however, it will not be possible to reach every point with positive probability. This is circumvented by choosing a measure which measures if a set is big enough for the walk to pass through that set with positive probability.

\begin{dfn}
    Let $\nu$ be a $\sigma$-finite, positive measure on $X$. We say that the chain is $\nu$\emph{-irreducible} if for any $x\in X$, for any $A\subset X$ with $\nu(A)>0$, there is some integer $n$ such that
    $$P^n(x,A)>0.$$
\end{dfn}

In the discrete case with $\nu$ being the counting measure, $\nu$-irreducibility amounts to usual irreducibility. This notion is rather weak and not canonical however, since in the discrete setting for example, irreducibility implies $\delta_x$ irreducibility for any $x$, and $\delta_x$-irreducibility is achieved as long as $x$ can be reached for any state, which is much weaker than actual irreducibility. We want something which also goes in the other direction: if $\nu(A)=0$, then the amount of starting points for which the walk visits $A$ is also small (we cannot ask this of every point since the walk starting at a point in $A$ very much visits $A$).

To remedy this, there is a notion of maximal irreducibility measure (class), making the above definition more canonical (see \cite[Proposition 4.2.2]{MT09}).

\begin{pro}\label{pro: maximal irreducibility measure}
    If $P$ is $\mu$-irreducible for some measure $\mu$, there is a measure $\nu$ such that 
    \begin{enumerate}[label=(\roman*)]
        \item $P$ is $\nu$-irreducible,
        \item If $\mu'$ is such that $P$ is $\mu'$-irreducible, then $\mu'\ll\nu$,
        \item If $\nu(A)=0$, then the set 
        $$\{x\in X,\m P_x(T_A<\infty)>0\}$$
        is $\nu$-negligible.
    \end{enumerate}
\end{pro}

Note that by the second point of Proposition \ref{pro: maximal irreducibility measure}, there is a \emph{unique} maximal measure class. 

\emph{From now on, when we say $P$ is $\nu$-irreducible, we always assume that $\nu$ is maximal in the sense of Proposition \ref{pro: maximal irreducibility measure}.}\medskip

When we defined balayage measures in the last section, we asked for the first hitting time to be finite a.s. in order to be able to define our objects.

\begin{dfn}\label{dfn: Harris recurrent set}
    A (measurable) set $A\subset X$ is said to be \emph{Harris-recurrent} if for any $x\in X$, 
    $$\m P_x(T_A<\infty)=1.$$
\end{dfn}

We would now like to have a way to tell when a subset is Harris-recurrent, and especially when big enough subsets are Harris-recurrent. This is where irreducibility comes in.

\begin{dfn}
    A $\nu$-irreducible chain $P$ is said to be recurrent in the sense of Harris, or  Harris-recurrent if for all $A$ such that $\nu(A)>0$, $A$ is Harris-recurrent.
\end{dfn}

For a discrete chain and $\nu$ the counting measure we recover the usual meaning of recurrence. Note also that since the maximal measure class is unique, the condition $\nu(A)>0$ does not depend on the particular choice of $\nu$.

However, this notion may be too strong for some purposes, because we require that it be reached with probability $1$ from \emph{every} point, but something bad could be happening on a small set. Hence a different notion of recurrence is the following:

\begin{dfn}
    A $\nu$-irreducible chain $P$ is said to be \emph{recurrent} if for all $A$ such that $\nu(A)>0$ and for all $x\in X$,
    $$\sum_{n=1}^\infty P^n(x,A)=\m E_x\left(\sum_{n=1}^\infty \m 1_{\{Z_n\in A\}}\right)=\infty.$$
\end{dfn}

Although it may seem unclear why the notions of recurrence and Harris recurrence could be different, they are in general, and the following Theorem both explains how that could be the case and explains why the two are not so different (see \cite[Theorem 9.1.5]{MT09}).

\begin{theo}\label{theo: recurrent=Harris+transient}
    Let $(Z_n)_n$ be a $\nu$-irreducible recurrent chain. Then there exists a unique decomposition 
    $$X=H\sqcup N,$$
    where $N$ is transient $\nu$-null, $H$ is absorbing (i.e. $P(x,H)=1$ for all $x\in H$) and $(Z_n)$ is Harris recurrent in restriction to $H$.
\end{theo}

Uniqueness is not stated in \cite{MT09}; however it is clear that if $H,H'$ are maximal absorbing sets which intersect on a positive measure set, they have to be equal. 

One way to understand this statement is that it could be that on some negligible set the walk escapes to infinity but still has a positive probability of going back to the reasonable set. This phenomenon can only happen if the walk is not homogeneous though.

\begin{cor}\label{cor: recur+transi->Harris}
    Let $(Z_n)_n$ be a $\nu$-irreducible recurrent chain, and assume that $\Aut(X,P)$ is acting transitively on $X$. Then $(Z_n)_n$ is Harris-recurrent.
\end{cor}

\begin{proof}
    The sets $H$ and $N$ in the decomposition of Theorem \ref{theo: recurrent=Harris+transient} are unique, and must thus be preserved by $\Aut(X,P)$, so that $gN=N$ for all $g\in G$. If $N$ were not empty, transitivity of $G$ implies that $N=GN=X$, contradicting the fact that $N$ is $\nu$-null.
\end{proof}

This Corollary will be useful for the proof of Theorem \ref{theointro: LS pour des marches sur des groupes lcsc}, when combined with the following criterion for recurrence (see \cite[Theorem 10.1.1]{MT09}:

\begin{pro}\label{pro: inv proba->recur}
    Assume that a probability measure $m$ is stationary, i.e. 
    $$\int_XP(x,A)\d m(x)=m(A)$$
    for all measurable $A$. Then $P$ is recurrent.
\end{pro}

The converse is also true, in the sense that for a recurrent $P$ there exists an invariant measure $m$ on $X$, although it may not be finite. It is unique up to rescaling however, and satisfies $m\succ \nu$ (see \cite[Theorem 10.0.1]{MT09}, \cite[Theorem 2.2.5]{Rev84}).

\subsection{Poisson-Furstenberg boundaries}

Theorem \ref{theointro: LS pour des marches sur des groupes lcsc} is about an isomorphism between \emph{Poisson-Furstenberg boundaries} of a group and a lattice of that group. We give here the two classical definitions in the context of a general Markov chain (see \cite{Kai95}).
\medskip

Let $(X,P)$ be a Markov chain as described before, let $o\in X$, and $(X^\m N,\m P_o)$ be the canonical realisation of that Markov chain. Let ${T:X^\m N\to X^\m N}$ be the time shift $T((x_n)_n)=(x_{n+1})_n$. Define $\Omega$ to be the space of ergodic components of $T$, which is naturally endowed with a projection $\pr :X^\m N\to\Omega$, and $\nu_o=\pr_*\m P_o$. 

The space $\Omega$ is called the \emph{Poisson-Furstenberg} boundary of $(X,P)$, and will often be denoted $B_{PF}(X,P)$ in what follows.\medskip

Another way to describe $\Omega$ is the following. Let $\mc I$ be the invariant $\sigma$-algebra, defined by the sets $A\in \mc F$ such that $T\inv (A)=A$. Since $X^\m N$ is standard Borel, Mackey's point realisation theorem gives the existence of a standard Borel space $(\Omega, \mc A)$ and a projection $\pr:X^\m N\to\Omega$ such that $\mc I=\pr^*(\mc A)$. \medskip

There is now a correspondence between bounded harmonic functions and $\Omega$. Indeed, for $h$ bounded harmonic, $(h(Z_n))_n$ converges almost surely by the martingale convergence theorem, and the limit must be $\mc I$-measurable and bounded. As such, it descends to a function in $\L^\infty(\Omega,\nu_o)$. Conversely, one can check that by irreducibility, all the $\nu_x$ are equivalent, so that $f\in \L^\infty(\Omega,\nu_o)$ is in all $\L^\infty(\Omega,\nu_x)$, and the map $h$ given by 
$$h(x)=\int_\Omega f\d\nu_x$$
can be checked to be $P$-harmonic and bounded. The two maps that have just been described are linear isomorphisms of the spaces 
$H^\infty(X,P)$ and $\L^\infty(\Omega, \nu_o)$ and are inverses of each other. They are moreover clearly $\Aut(X,P)$-equivariant.

Hence in order to prove isomorphisms between the Poisson boundaries of two random walks, one may try to find an isomorphism between their spaces of bounded harmonic functions instead.

\section{A discretization method}\label{sec: discretisation}
In this section, we prove Theorem \ref{theointro: LS pour des marches sur des groupes lcsc} through a discretization procedure, as well as auxiliary results. We follow the paper of Lyons and Sullivan \cite{LS84}, although many different authors have provided slightly different descriptions of the same kind of ideas, for instance \cite{Fur63,Anc90,Kai92,BL96,KM96} to name a few.
\subsection{A discrete-time Lyons-Sullivan procedure}
We start by the following definition, which will be our way of getting a \emph{Harnack inequality} in this context.

\begin{dfn}\label{dfn: condition (H)}
    Let $G$ be a locally compact, second countable group. Let $\mu$ be a probability measure on $G$.
    
    We say $\mu$ \emph{satisfies condition (H)} if there exist a relatively compact neighbourhood $U$ of $1$ and a constant $C$ with $0<C<1$ such that for all $g\in U$, 
    $$\frac{\d (g_*\mu)}{\d\mu}(x)\ge C$$
    for $\mu$-a.e. $x\in G$.
\end{dfn}

This condition is quite strong; it may not hold even for compactly supported measures. We will go over methods for weakening this hypothesis in Section \ref{ssec: Weaker H}.\medskip

\begin{rem}
    Note that $U$ may be shrunk to be made symmetric. In this case, for any $g\in U$,
    $$\frac{\d g\inv_*\mu}{\d\mu}(x)=\frac{\d \mu}{\d g_*\mu}(gx)\ge C$$
    for $\mu$-a.e. $x\in G$. This gives an upper bound for the Radon-Nikodym derivatives as well, which will be useful later on in Section \ref{ssec: stopping times everywhere}.
\end{rem}

This is the actual way condition (H) will be used in what follows.
\begin{lem}\label{lem: convo mu -> RN cool}
    Let $\mu$ be a measure satisfying (H). For any nonnegative finite measure $\nu$ supported on $U$, 
    $$\frac{\d(\nu*\mu)}{\d\mu}(x)\ge C\nu(G)$$
    for $\mu$-a.e. $x\in G$. 
\end{lem}

\begin{proof}
    For $\nu=\delta_g$ with $g\in U$, it is simply condition (H).
    
    For the general case, observe that 
    $$\frac{\d(\nu*\mu)}{\d\mu}(x)=\int_G\frac{\d(\delta_g*\mu)}{\d\mu}(x)\d\nu(g)\ge \int_GC\d\nu(g)=C\nu(G)$$
    as claimed.
\end{proof}

The following is Theorem \ref{theointro: LS pour des marches sur des groupes lcsc}.

\begin{theo}\label{theo: BPF lattice in lcsc group}
    Let $G$ be a locally compact, second countable group. Let $\mu$ be an admissible probability measure on $G$ which satisfies (H).
    Let $\Gamma<G$ be a lattice. Then there is an admissible measure $\nu$ on $\Gamma$ such that 
    $$B_{PF}(G,\mu)\simeq B_{PF}(\Gamma,\nu).$$
    Moreover, the isomorphism is $\Gamma$-equivariant, and $\nu$ can be chosen to be symmetric whenever $\mu$ is.
\end{theo}

Note first that since $G$ has a lattice, it is unimodular \cite{Rag72}, so that we will make no difference between the left and right Haar measure.

Let $(Z_n)_{n\in\m N}$ be the Markov chain associated to $\mu$ as in Example \ref{ex: Markov kernels} \ref{ex: Markov kernels; on groups}, i.e. a sequence of random variables in $G$ such that $Z_n\inv Z_{n+1}$ are all i.i.d of law $\mu$. 

We are interested in proving that for some (in fact any) neighbourhood $U$ of $1$ in $G$, the set $\Gamma U$ is recurrent for the walk. For this, we will prove that the set of \emph{right} $\Gamma$ classes is recurrent in the quotient, hence the need for right quotients for the following part.

Let $m$ be the $G$-invariant probability measure on $G/\Gamma$, and let $\pi:G\to G/\Gamma$ be the canonical projection.

\begin{pro}\label{pro: lattice in group is Harris}
     For any Borel set $A\subset G/\Gamma$ with $m(A)>0$ and for every $x\in G/\Gamma$, 
    $$\m P_x(\exists n,\pi(Z_n)\in A)=1.$$
\end{pro}

\begin{proof}
    Let $P$ be the Markov kernel associated to $\pi(Z_n)$, i.e. 
    $$P(x,A)=(\mu*\delta_x)(A).$$
    Note that since $\mu$ is admissible, the chain $P$ is $m$-irreducible, and $m$ is in fact invariant for $P$. Indeed, by Fubini, for any Borel subset $A$ of $G/\Gamma$,
    \begin{eqnarray*}
        \int_{G/\Gamma}P(x,A)\d m(x)
        &=&\int_{G/\Gamma}(\mu*\delta_x)(A)\d m(x)\\
        &=&\int_{G/\Gamma}\int_G g_*\delta_x(A)\d\mu(g)\d m(x)\\
        &=&\int_G m(A g\inv)\d\mu(g)=m(A).
    \end{eqnarray*}
    Since $m$ is of finite volume, by Proposition \ref{pro: inv proba->recur} and Corollary \ref{cor: recur+transi->Harris}, $P$ is Harris-recurrent, which directly implies the statement.
\end{proof}
Recall the definition of the first hitting time of a Borel set $A\subset X$, defined by
$$T^A=\inf\{n\in \m N, Z_n\in A\},$$
as well as the balayage measure $\nu^A$ for Harris-recurrent $A$, which is the distribution of $Z^A:=Z_{T^A}$ with starting distribution $\nu$, i.e.
$$\nu^A(B)=\m P_\nu(Z_{T^A}\in B).$$

Fix for the rest of the proof $U$ and $C$ given by condition (H) (see Definition \ref{dfn: condition (H)}). Up to shrinking $U$, assume moreover that $\gamma U\cap U=\emptyset$ for any $\gamma\in\Gamma$ with $\gamma\ne 1$, which is possible since $\Gamma$ is a discrete subgroup, and that $U$ is symmetric.

\begin{proof}[Proof of Theorem \ref{theo: BPF lattice in lcsc group}]
    When $G$ is discrete, simply letting $\nu=\delta_e^\Gamma$ works; see for instance \cite{KW02,Der25}. We will hence assume from now on that $G$ is not discrete, and hence that the Haar measure of $G$ is non atomic. This will be useful later on.\medskip
    
    The idea of the proof is the same as that of Lyons-Sullivan, in that we will construct two sequences of measures inductively, and transfer more and more of the mass from the first sequence to the second. More precisely, given a finite positive measure $\nu_n$ on $G$, we construct two measures $\nu_{n+1}$, $\tau_{n+1}$ on $G$, such that 
    \begin{enumerate}
        \item $\tau_{n+1}$ is supported on $\Gamma$,
        \item $\nu_n(h)=\nu_{n+1}(h)+\tau_{n+1}(h),$
        \item $\tau_{n+1}(G)=C\nu_n(G)$
    \end{enumerate}
    for any bounded positive $\mu$-harmonic function $h$.\medskip

    Start by moving all mass away from $\Gamma$, by letting 
    $$\nu_n'=\nu_n+\sum_{\gamma\in\Gamma}\nu_n(\gamma)(\gamma_*\mu-\delta_\gamma).$$
    Note that this does not change the value of $\nu_n(h)=\nu_n'(h)$ for all bounded positive harmonic $\mu$-harmonic function $h$. Now since $\Gamma U$ is Harris-recurrent by Proposition \ref{pro: lattice in group is Harris}, let $\tilde\nu=(\nu_n')^{\Gamma U}$, which is a well defined measure supported on $\Gamma U$, and for which $\nu_n(h)=\nu_n'(h)=\tilde \nu(h)$ for any bounded $\mu$-harmonic function $h$ by Lemma \ref{lem: properties of balayage}. For $\gamma\in\Gamma$, let $\rho_\gamma$ be the restriction of $\tilde\nu$ to $\gamma U$, that is, $\rho_\gamma(A)=\tilde\nu(A\cap \gamma U)$. Let also $m_\gamma=\rho_\gamma(G)=\tilde\nu(\gamma U)$. Note that since the $\gamma U$ are disjoint, 
    $$\tilde\nu=\sum_{\gamma\in\Gamma}\rho_\gamma.$$
    Moreover, since $\rho_\gamma$ is supported on $\gamma U$, the measure $(\gamma\inv)_*\rho_\gamma$ is supported on $U$ and
    \begin{equation}
    \frac{\d(\rho_\gamma*\mu)}{\d(\delta_\gamma*\mu)}(x)=\frac{\d((\gamma\inv)_*\rho_\gamma*\mu)}{\d\mu}(\gamma\inv x)\ge Cm_\gamma\label{eqn: control on RN on gamma U}
    \end{equation}
    for $\gamma_*\mu$-a.e. $x\in G$ by Lemma \ref{lem: convo mu -> RN cool}.\medskip

    Note that for any bounded $\mu$-harmonic function $h$, 
    \begin{eqnarray*}
        \nu_n(h)&=&\sum_{\gamma\in\Gamma}(\rho_\gamma*\mu)(h)\\
        &=&\sum_{\gamma\in\Gamma}\left[\left(\rho_\gamma*\mu-Cm_\gamma\delta_\gamma*\mu\right)(h)+Cm_\gamma h(\gamma)\right].
    \end{eqnarray*}
    Hence if we let 
    $$\nu_{n+1}=\sum_{\gamma\in\Gamma}\left(\rho_\gamma*\mu-Cm_\gamma\delta_\gamma*\mu\right)\quad\text{and}\quad \tau_{n+1}=\sum_{\gamma\in\Gamma}Cm_\gamma\delta_\gamma,$$
    then by (\ref{eqn: control on RN on gamma U}), $\nu_{n+1}$ is a positive measure, and we have $\nu_n(h)=\nu_{n+1}(h)+\tau_{n+1}(h)$ for any bounded $\mu$-harmonic function $h$.

    Finally, since 
    $$\sum_{\gamma\in\Gamma}m_\gamma=\sum_{\gamma\in\Gamma}\rho_\gamma(G)= \tilde \nu(G)=\nu_n(G),$$
    we have $\tau_{n+1}(G)=C\nu_n(G)$ as claimed.\medskip

    Fix now a finite positive measure $\nu=\nu_0$ on $G$. If $h$ is a bounded positive $\mu$-harmonic function, then
    $$\nu_0(h)=\nu_n(h)+\sum_{k=1}^n\tau_k(h).$$
    Since $\nu_n(G)=(1-C)^n\nu_0(G)$, for any bounded harmonic function $h$ we have $|\nu_n(h)|\le (1-C)^n\|h\|_\infty\nu_0(G)\tendv 0$. Hence, the measure $\Phi(\nu_0):=\sum_{k\ge 1}\tau_k$ is finite and such that $\Phi(\nu_0)(h)=\nu_0(h)$ for any bounded $\mu$-harmonic measure on $G$.
    
    \begin{lem}\label{lem: propriétés de Phi}
        The map $\Phi:\nu_0\mapsto\sum\tau_k$ can be extended by $\Phi(\nu)=\Phi(\nu^+)-\Phi(\nu^-)$ to any signed measure. It is then linear, $\Gamma$-equivariant, continuous for the norm of total variation, and satisfies $\Phi(\nu_0)=\Phi(\nu_1)+\tau_1$.
    \end{lem}

    \begin{proof}
        It is straightforward to check equivariance and that for positive measures $\nu,\nu'$ and some $\alpha>0$, $\Phi(\nu+\nu')=\Phi(\nu)+\Phi(\nu')$ and $\Phi(\alpha\nu)=\alpha\Phi(\nu)$. The extension of $\Phi$ to signed measures is then clearly linear. Moreover, $\Phi(\nu)(G)=\nu(G)$, so that $\Phi$ is $1$-Lipschitz over the cone of positive measures for the norm of total variation, and hence over all signed measures. Finally, the formula $\Phi(\nu_0)=\Phi(\nu_1)+\tau_1$ comes from the definition of $\Phi$ and the fact that the definition of $\nu_{n+1}$ only depends on $\nu_n$ and not any of the previous $\nu_i$'s.
    \end{proof}

    We are now able to construct measures supported on $\Gamma$ which have the same Poisson boundary as the walk on $G$. For $x\in G$, letting $\nu_0=\delta_x$ lets us define $p_x$ by $p_x=\Phi(\delta_x)$. The measure $p_x$ is a probability measure supported on $\Gamma$, so we can define a random walk $\hat P$ on $\Gamma$ by letting $\hat p(x,y)=p_x(y)$.

    We now do the computations to check that the Poisson boundaries actually coincide. We start by the following definition, which will aid us in the computations.\medskip

    For $f\in \ell^\infty(\Gamma),$ define $h_f:G\to\m R$ by $h_f(x)=p_x(f)$.

    \begin{lem}\label{lem: Phi et harmonicité}
        For $\theta\in\rm{Prob}(G)$ we have
        $$\theta(h_f)=\Phi(\theta)(f).$$
    \end{lem}

    \begin{proof}
        Since $\Phi$ is linear continuous,
        $$\int_Gp_y\d\theta(y)=\int_G\Phi(\delta_y)\d\theta(y)= \Phi\left(\int_G\delta_y\d\theta(y)\right)=\Phi(\theta),$$
        so that
        $$\theta(h_f)=\int_Gh_f(y)\d\theta(y)= \int_Gp_y(f)\d\theta(y)= \Phi(\theta)(f).$$
    \end{proof}

    \begin{lem}
        For a bounded positive $\hat P$-harmonic map $f$, the map $h_f$ is (bounded positive) $\mu$-harmonic.
    \end{lem}
    
    \begin{proof}
        We start by computing the quantity $\Phi(x_*\mu)(f)$. For this we consider the following three cases:
    
        \begin{itemize}
            \item If $x=\gamma_0\in\Gamma$, then the measure $(\gamma_0)_*\mu$ does not charge $\Gamma$ since $\mu$ is absolutely continuous with respect to Haar, so for both $\nu_0=\delta_{\gamma_0}$ and $\nu_0=(\gamma_0)_*\mu$, $\nu_0'=(\gamma_0)_*\mu$, so that 
            $$\Phi((\gamma_0)_*\mu)=\Phi(\delta_{\gamma_0})=p_{\gamma_0}.$$
            But then by $\hat P$-harmonicity of $f$, 
            $$\Phi((\gamma_0)_*\mu)(f)=p_{\gamma_0}(f)=f(\gamma_0).$$
            
            \item If $x\in\Gamma U$, say $x\in\gamma_0 U$ but $x\not\in\Gamma$, then for $\nu_0=\delta_x$, we have  $(\nu_0')^{\Gamma U}=\delta_x^{\Gamma U}=\delta_x,$ $\rho_\gamma=\m 1_{\gamma_0}(\gamma)\delta_x$ and $m_\gamma=\m 1_{\gamma_0}(\gamma)$. 

            In particular, $\nu_1=x_*\mu-C(\gamma_0)_*\mu$ and $\tau_1=C\delta_{\gamma_0}$. By Lemma \ref{lem: propriétés de Phi}, $\Phi(\nu_1)=\Phi(\nu_0)-\tau_1$, which gives
            \begin{eqnarray*}
                \Phi(x_*\mu)(f)&=& \Phi(\delta_x)(f)-C\delta_{\gamma_0}(f)+C\Phi((\gamma_0)_*\mu)(f)\\
                &=&p_x(f)-Cf(\gamma_0)+Cf(\gamma_0)=p_x(f).
            \end{eqnarray*}
            \item If $x\not\in \Gamma U$, then $\delta_x^{\Gamma U}=(x_*\mu)^{\Gamma U}$ since $T^{\Gamma U}\ge 1$ $\m P_x$-a.s., so that 
            $$p_x=\Phi(\delta_x)=\Phi(x_*\mu).$$
        \end{itemize}

        In all cases, we have $\Phi(x_*\mu)(f)=h_f(x)$, and we are now able to conclude by Lemma \ref{lem: Phi et harmonicité}:
        $$(x_*\mu)(h_f)=\Phi(x_*\mu)(f)=h_f(x),$$
        which is precisely $\mu$-harmonicity of $h_f$.
    \end{proof}
    
    The fact that the restriction of $h_f$ to $\Gamma$ coincides with $f$ comes from the fact that $h_f(\gamma)=p_\gamma(f)=f(\gamma)$ since $f$ is $\hat P$-harmonic. Moreover, for any $\mu$-harmonic $h$, 
    $$p_x(h|_\Gamma)=p_x(h)=\Phi(\delta_x)(h)=\delta_x(h)=h(x).$$ 
    The two maps we have defined are thus inverses of each other, and are clearly $\Gamma$-equivariant.\medskip

    Finally, $p_1$ is admissible. Indeed, let 
    $$Y=\bigcup_{n\ge 0}\rm{supp}(p_1^{*n}).$$
    If $Y\ne\Gamma$, let $\gamma\in\Gamma-Y$. Note that $\m 1_Y$ is $p_1$-harmonic. In particular, by what we have proven before, $h:x\mapsto p_x(\m 1_Y)$ is $\mu$-harmonic. Recall also that for $x\in U\setminus\{e\}$ and $\nu_0=\delta_x$, we had $\tau_1=C\delta_e$, so that 
    $$h(x)=p_x(\m 1_{Y})\ge p_x(e)\ge\tau_1(e)=C.$$
    But if we now take some $n$ such that $\mu^{*n}(\gamma\inv (U\setminus\{e\}))>0$ (which exists because $\mu$ is admissible), we get 
    $$0=\m 1_{Y}(\gamma)=h(\gamma)=\delta_\gamma*\mu^{*n}(h)\ge C\mu^{*n}(\gamma\inv (U\setminus\{e\}))>0$$
    which is a contradiction, and $p_1$ is thus admissible. \bigskip

    Since $\Phi$ is $\Gamma$-equivariant, $p_\gamma=\gamma_*p_1$, and the random walk $\hat P$ is simply the walk associated to $p_1$. Hence, the measure $\nu=p_1$ on $\Gamma$ does verify 
    $$H^\infty(G,\mu)\simeq H^\infty(\Gamma,\nu)$$
    and 
    $$B_{PF}(G,\mu)\simeq B_{PF}(\Gamma,\nu)$$
    where the isomorphisms are $\Gamma$-equivariant.
\end{proof}

\subsection{Description of the walk on $\Gamma$ by stopping times}\label{ssec: stopping times everywhere}
Informally, the discretization process can be thought of in the following way: start the walk with some distribution $\nu_0$. Walk until you hit $\Gamma U$, with the slight subtlety that if we started at some point in $\Gamma$ we cannot stop right away. Then the fact that we do not just have $\Phi(\nu_0)=\sum_\gamma m_\gamma \delta_\gamma$ can be thought of as there only being a given probability that one would stop the first time hitting $\Gamma U$, and if this is not achieved then start walking again until you hit $\Gamma U$ for the second time, and if the coin flip fails again, try at the third. The probability of stopping at any given hitting time is not very simple however; it needs to simulate the formula 
$$\nu_1=\tilde \nu*\mu-C\sum_\gamma m_\gamma\gamma_*\mu.$$
In order to do this, the probability to stop at some hitting time has to be skewed by some Radon-Nikodym derivative. The actual (quite technical) definition goes as follows, with the notations made to mirror those found in \cite{Kai92}.

Start by defining two sequences of stopping times $(R_n)_{n\ge 1}$ and $(S_n)_{n\ge 0}$ inductively for any $\omega\in G^\m N$ by 
$$S_0(\omega)=\m 1_{\{\omega_0\in\Gamma\}},$$
and then 
$$R_n(\omega)=\inf\{k\ge S_{n-1}(\omega),\omega_k\in\Gamma U\}$$
and 
$$S_n(\omega)=R_n(\omega)+1.$$
The reason why $S_n$ is deterministic given $R_n$ in our setup (and not in the original Lyons-Sullivan setup) is that we only need to walk the measure one step further, instead of doing another balayage away from $\Gamma U$.

Note that all the $S_n$ are finite a.s. by Proposition \ref{pro: lattice in group is Harris}.

One way to model the skewed coin flips to determine at which hitting time we stop is to extend the probability space $(G^\m N,\m P_\nu)$ to the space $(\tilde\Omega,\tilde{\m P}_\nu)=(G^\m N\times [0,1]^\m N,\m P_\nu\otimes \rm{Leb}^{\otimes \m N})$. If we let $\gamma_k(\omega)$ be the only $\gamma\in\Gamma$ such that $\omega_{R_k(\omega)}\in\gamma U$ and $\rho_{\gamma,k}$ be the distribution of $Z_{R_k}$ restricted to $\gamma U$ (that is, $\rho_{\gamma,k}(A)=\m P(Z_{R_k}\in A\cap \gamma U)$), we can then define the stopping time $N$ by 
$$N((\omega_n)_n,(\alpha_n)_n)=\inf\left\{k\ge 1, C\|\rho_{\gamma_k(\omega),k}\|\frac{\d(\gamma_k(\omega))_*\mu}{\d \rho_{\gamma_k(\omega),k}*\mu}(\omega_{S_k(\omega)})>\alpha_k\right\}.$$

Note that since $U$ is assumed to be symmetric, for any $\gamma\in\Gamma$,
$$C\|\rho_{\gamma,k}\|\le \frac{\d \rho_{\gamma,k}*\mu}{\d\gamma_*\mu}(x)\le \frac1C\|\rho_{\gamma,k}\|,$$
so that 
$$C^2\le C\|\rho_{\gamma_k(\omega),k}\| \frac{\d(\gamma_k(\omega))_*\mu}{\d \rho_{\gamma_k(\omega),k}*\mu}(\omega_{S_k(\omega)})\le 1$$
$\m P$-a.e. In particular, $N\le \tilde N$ a.s., where
$$\tilde N((x_n)_n,(\alpha_n)_n)=\inf\left\{k, C^2>\alpha_k\right\}.$$
Note that the law of $\tilde N$ for $\tilde{\m P}_\eta$ is geometric of parameter $C^2$. In particular, $N$ is finite a.s. and has an exponential moment.

\begin{lem}\label{lem: distri de ZRN est nu}
    For any finite measure $\nu_0$, we have
    $$\Phi(\nu_0)(\gamma)=\tilde{\m P}_{\nu_0}(Z_{R_N}\in\gamma U).$$
    More precisely, for any bounded measurable $f$ and any $\gamma\in \Gamma$,
    $$\tilde{\m E}_{\nu_0}(\m 1_{\{Z_{R_N}\in\gamma U\}}f(Z_{{S_N}}))=\Phi(\nu_0)(\gamma)\times (\gamma_*\mu)(f).$$
\end{lem}

\begin{proof}
    Note that the first statement follows from the second one by taking $f\equiv 1$.

    Now
    \begin{eqnarray*}
    &&\tilde{\m E}_{\nu_0}(f(Z_{S_1})\m 1_{\{N=1\}}\m 1_{\{Z_{R_1}\in\gamma U\}})\\
    &&=\int_{G^\m N\times [0,1]^\m N}f(Z_{S_1}(\omega))\m 1_{\{N(\omega,\alpha)=1\}}\m 1_{\{Z_{R_1}\in \gamma U\}}\d\tilde{\m P}_{\nu_0}(\omega,\alpha)\\
    &&=\int_{G^\m N}f(Z_{S_1}(\omega))\m 1_{\{Z_{R_1}\in \gamma U\}}\int_{[0,1]^\m N}\m 1_{\{N(\omega,\alpha)=1\}}\d\rm{Leb}^{\otimes\m N}(\alpha)\d\m P_{\nu_0}(\omega)\\
    &&=\int_{G^\m N}f(Z_{S_1}(x)) \m 1_{\{Z_{R_1}\in \gamma U\}}C\|\rho_{\gamma,1}\|\frac{\d\gamma_*\mu}{\d \rho_{\gamma,1}*\mu}(Z_{S_1}(\omega))\d\m P_{\nu_0}(\omega)\\
    &&=C\|\rho_{\gamma,1}\|\int_G\int_Gf(xy)\frac{\d\gamma_*\mu}{\d\rho_{\gamma,1}*\mu}(xy)\d\rho_{\gamma,1}(x)\d\mu(y)\\
    &&=C\|\rho_{\gamma,1}\|\int_Gf(x) \frac{\d\gamma_*\mu}{\d\rho_{\gamma,1}*\mu}(x)\d\rho_{\gamma,1}*\mu(x)\\
    &&=C\|\rho_{\gamma,1}\|\gamma_*\mu(f)=\tau_1(\gamma)\gamma_*\mu(f).
\end{eqnarray*}

Note that $\tilde{\m E}_{\nu_0}(f(Z_{S_1})\m 1_{\{Z_{R_1}\in\gamma U\}})=\rho_{\gamma,1}*\mu(f)$, so summing over $\Gamma$ we also obtain
$$\nu_1(f)=\tilde{\m E}_{\nu_0}(f(Z_{S_1})\m 1_{\{N>1\}}).$$
Let $\theta$ be the shift on $G^\m N$: 
$$\theta((\omega_0,\omega_1,\ldots))=(\omega_1,\omega_2,\ldots).$$
Clearly, $Z_n\circ \theta^p=Z_{n+p}$. For $T:\Omega\to\m N$ an a.s. finite random variable, let $\theta_T$ denote the random variable $\theta_T(\omega)=\theta^{T(\omega)}(\omega)$.

Note that $R_{n+1}=S_1+R_n\circ\theta_{S_1}$ so that $Z_{R_{n+1}}=Z_{R_n}\circ\theta_{S_1}$, and likewise $Z_{S_{n+1}}=Z_{S_n}\circ\theta_{S_1}$. Moreover,
$$\{N=n+1\}=\{N\circ\theta_{S_1}=n\}\cap\{N>n\}.$$ 
As such by the strong Markov property,
\begin{eqnarray*}
    &&\tilde{\m E}_{\nu_0}(f(Z_{S_2})\m 1_{\{N=2\}}\m 1_{\{Z_{R_2}\in\gamma U\}})\\
    &&=\tilde{\m E}_{\nu_0}\Big(\tilde{\m E}_{\nu_0}(f(Z_{S_2})\m 1_{\{N=2\}}\m 1_{\{Z_{R_2}\in\gamma U\}}\mid \mc F_{S_1})\Big)\\
    &&=\tilde{\m E}_{\nu_0}\Big(\m 1_{\{N>1\}}\tilde{\m E}_{\nu_0}(f(Z_{S_1}\circ\theta_{S_1})\m 1_{\{N\circ\theta_{S_1}=1\}}\m 1_{\{Z_{R_1}\circ\theta_{S_1}\in\gamma U\}}\mid \mc F_{S_1})\Big)\\
    &&=\tilde{\m E}_{\nu_0}\Big(\m 1_{\{N>1\}}\tilde{\m E}_{Z_{S_1}}(f(Z_{S_1})\m 1_{\{N=1\}}\m 1_{\{Z_{R_1}\in\gamma U\}})\Big)\\
    &&=\tilde{\m E}_{\nu_1}(f(Z_{S_1})\m 1_{\{N=1\}}\m 1_{\{Z_{R_1}\in\gamma U\}})\\
    &&=\tau_2(\gamma)\gamma_*\mu(f)
\end{eqnarray*}
where the second to last equality comes from the formula $\nu_1(f)=\tilde{\m E}_{\nu_0}(f(Z_{S_1})\m 1_{\{N>1\}})$ applied to $f$ being the conditional expectation. 

The exact same argument works for $k\ge 3$, so by induction it also holds that 
$$\tilde{\m E}_{\nu_0}(f(Z_{S_k})\m 1_{\{N=k\}}\m 1_{\{Z_{R_k}\in\gamma U\}})=\tau_k(\gamma)\gamma_*\mu(f)$$
for all $k\ge 1$. Summing over all $k$, we obtain 
$$\tilde{\m E}_{\nu_0}(\m 1_{\{Z_{R_N}\in\gamma U\}}f(Z_{{S_N}}))=\Phi(\nu_0)(\gamma)\times (\gamma_*\mu)(f),$$
as hoped.
\end{proof}

We are now interested in seeing if the random walk on $\Gamma$ can be modelled on that of $G$. The answers turns out to be yes, and to that end we not only record the first moment the coin flip hits, but the $k$-th moment that this happens; this amounts to letting inductively $N_1=N$, and 
$$N_{k+1}((\omega_n)_n,(\alpha_n)_n)=\inf\left\{l> N_k, C\|\rho_{\gamma_l(\omega),l}\|\frac{\d(\gamma_l(\omega))_*\mu}{\d \rho_{\gamma_l(\omega),l}*\mu}(\omega_{S_l(\omega)})>\alpha_l\right\}.$$

\begin{pro}\label{pro: approximation by N_k}
    For any $n\in \m N$ and $\gamma_1,\ldots,\gamma_n\in\Gamma$,
    $$\tilde{\m P}_e(Z_{R_{N_1}}\in \gamma_1U,\ldots, Z_{R_{N_n}}\in\gamma_nU)=\m Q_e(Y_1=\gamma_1,\ldots,Y_n=\gamma_n),$$
    where $\m Q_e$ is the canonical probability measure on $\Gamma^\m N$ associated to $\nu$ and $Y_i:\Gamma^\m N\to\Gamma$ is the $i$-th coordinate projection (see Definition \ref{prodfn: canonical Markov chain}).
\end{pro}

\begin{proof}
    Much like for the proof of Lemma \ref{lem: distri de ZRN est nu} we prove by induction that for all $n$ and all $f$, 
    $$\tilde{\m E}_e(\m 1_{\{Z_{R_{N_1}}\in \gamma_1U,\ldots, Z_{R_{N_n}}\in\gamma_nU\}}f(Z_{S_{N_n}}))=\m Q_e(Y_1=\gamma_1,\ldots,Y_n=\gamma_n)(\gamma_n)_*\mu(f).$$ 
    
    For $n=1$, 
    $$\tilde{\m E}_e(\m 1_{\{Z_{R_N}\in \gamma U\}}f(Z_{S_N}))=\nu(\gamma)\gamma_*\mu(f)$$
    is Lemma \ref{lem: distri de ZRN est nu}. Assume now that it is true for some $n\ge1$. For simplicity, for $k\ge 1$ let $a_k$ denote $\m Q_e(Y_1=\gamma_1,\ldots,Y_k=\gamma_k)$. Then
    \begin{eqnarray*}
        (*)&=&\tilde{\m E}_e(\m 1_{\{Z_{R_{N_1}}\in \gamma_1U,\ldots, Z_{R_{N_{n+1}}}\in\gamma_{n+1}U\}}f(Z_{S_{N_{n+1}}}))\\
        &=&\tilde{\m E}_e(\tilde{\m E}_e(\m 1_{\{Z_{R_{N_1}}\in \gamma_1U,\ldots, Z_{R_{N_{n+1}}}\in\gamma_{n+1}U\}}f(Z_{S_{N_{n+1}}})\mid \mc F_{S_{N_n}}))\\
        &=&\tilde{\m E}_e(\m 1_{\{Z_{T_{N_1}}\in \gamma_1U,\ldots, Z_{T_{N_n}}\in\gamma_nU\}}\tilde{\m E}_{Z_{T_{N_n}+1}}(\m 1_{\{Z_{T_N}\in \gamma_{n+1} U\}}f(Z_{T_N+1})))
    \end{eqnarray*}
    by the strong Markov property. Now we may apply the induction hypothesis with $f$ being the conditional expectation above, to get that
    \begin{eqnarray*}
        (*)&=&a_n\tilde{\m E}_{(\gamma_n)_*\mu}(\m 1_{\{Z_{T_N}\in \gamma_{n+1} U\}}f(Z_{T_N+1}))\\
        &=&a_n\Phi((\gamma_n)_*\mu)(\gamma_{n+1})\times(\gamma_{n+1})_*\mu(f)\\
        &=&a_n\Phi(\delta_{\gamma_n})(\gamma_{n+1})\times(\gamma_{n+1})_*\mu(f)\\
        &=&a_{n+1}(\gamma_{n+1})_*\mu(f)
    \end{eqnarray*}
    as hoped.
\end{proof}

\subsection{Weakening condition (H)}\label{ssec: Weaker H}

As noted after Definition \ref{dfn: condition (H)}, condition (H) is quite strong, and many classes of measures to which we would like to apply Theorem \ref{theo: BPF lattice in lcsc group} do not satisfy it, including most compactly supported measures. One solution is to modify the measure slightly, so that it keeps the same Poisson boundary but now satisfies condition (H). 

Fix some admissible $\mu\in\Prob(G)$, and let 
$$\hat\mu:=\sum_{n\ge 1}2^{-n}\mu^{*n}\in\Prob(G).$$
The measure $\hat\mu$ is an example of a \emph{sampling measure} (see for instance \cite[Section 5.5.1]{MT09}). It satisfies moreover the second goal of the modified measure we mentioned above:

\begin{pro}\label{pro: mu chapeau conserve la frontière de Poisson}
    We have the following $G$-equivariant isomorphism:
    $$B_{PF}(G,\mu)\simeq B_{PF}(G,\hat\mu).$$
\end{pro}

\begin{proof}
    This is an application of \cite[Lemma 6.3]{FT22} with $\alpha=\beta=\mu/2$.
\end{proof}

\begin{rem}
    It may also be that in specific cases, $\hat\mu$ may not satisfy condition (H), but another measure $\mu$ does; \cite[Proposition 4.7]{BHO23} or \cite[Lemma 6.3]{FT22} could then come in handy.
\end{rem}

We now give sufficient criteria for $\hat\mu$ to satisfy condition (H).

\begin{pro}\label{pro: condition pour que (H) ok}
    If there exist some integers $k\ge 0$, some constant $C'\in]0,1[$ and an open neighbourhood $U$ of $1$ such that for all $g\in U$, for $\mu$-a.e. $x\in G$,
    $$\frac{\d g_*\mu^{*k+1}}{\d\mu}(x)\ge C',$$
    then $\hat\mu$ satisfies condition (H).
\end{pro}

\begin{proof}
    Let $A\subset G$ be Borel. Then for any $g\in U$ and $l\ge 1$,
    \begin{eqnarray*}
        \mu^{*k+l}(g\inv A)&=&\int_G\mu^{*k+1}(g\inv A h\inv)\d\mu^{*l-1}(h)\\
        &\ge&C'\int_G\mu(Ah\inv)\d\mu^{*l-1}(h)=C'\mu^{*l}(A),
    \end{eqnarray*}
    hence for all $l\ge 1$, 
    $$\frac{\d g_*\mu^{*k+l}}{\d\mu^{*l}}(x)\ge C'.$$
    Now for $g\in U$ and $A\subset G$ Borel,
    \begin{eqnarray*}
        \hat\mu(g\inv A)&=&\sum_{n\ge 1}2^{-n}\mu^{*n}(g\inv A)\\
        &\ge&\sum_{n\ge k+1}2^{-n}\mu^{*n}(g\inv A)\\
        &\ge&\sum_{n\ge k+1}2^{-n}C'\mu^{*n-k}(A)\\
        &=&\sum_{n\ge 1}2^{-n-k}C'\mu^{*n}(A)=2^{-k}C'\hat\mu(A)
    \end{eqnarray*}
    Hence $\hat\mu$ satisfies (A) with $C=2^{-k}C'$ and the same $U$.
\end{proof}

\begin{cor}\label{cor: supp cpt-> mu chapeau a (H)}
    If $\lambda$ denotes the Haar measure on $G$, if $\mu\ll\lambda$ and $f=\d\mu/\d\lambda$ is continuous with compact support, $\hat\mu$ satisfies condition (H).
\end{cor}

\begin{proof}
    Since $\mu$ is generating, for some compact neighbourhood $U$ of $1$, we can find some $n\in\m N$ such that $\supp f^{*n}\supset U\inv\supp f$. Since $f^{*n}$ is continuous and positive on $U\inv\supp f$ there is an $\epsilon>0$ such that for any $y\in U\inv\supp f$, $f^{*n}(y)\ge\epsilon$ by compactness. Similarly, there is some $C''$ positive finite such that $f(x)\le C''$. But now,
    $$\frac{\d g_*\mu^{*n}}{\d\mu}(x)=\frac{f^{*n}(g\inv x)}{f(x)}\ge \frac\epsilon{C''}>0,$$
    and we are done by applying Proposition \ref{pro: condition pour que (H) ok}.
\end{proof}

\begin{rem}
    Some additional transformations of measures preserving the Poisson boundary are described in \cite{FT22}. In particular, when starting with a spread-out measure, one can assume without loss of generality that that measure is absolutely continuous with respect to Haar, and that its Radon-Nikodym derivative is bounded. Up to convolving that measure with itself once more, one can in fact assume that the Radon-Nikodym derivative is continuous. This is not enough to be able to apply Corollary \ref{cor: supp cpt-> mu chapeau a (H)} however, as the measure they construct will always have full support, and it is not clear if one can always obtain a measure that has (H) from any spread-out measure without changing the Poisson boundary.
\end{rem}

We end this section by recording that (H) is stable under products.

\begin{lem}\label{lem: stabilité de (H) sous produit}
    Let $G,G'$ be lcsc groups, and $\mu,\mu'$ be measures on $G$ and $G'$ respectively. If $\mu$ and $\mu'$ both satisfy (H), then so does $\mu\otimes \mu'$ on $G\times G'$.
\end{lem}

\begin{proof}
    Let $U,C,U',C'$ be the open sets and the constants given by (H) for $\mu$ and $\mu'$.
    
    For $(g,g')\in U\times U'$ and $(x,x')\in G\times G'$, we have
    $$\frac{\d(g,g')_*(\mu\otimes\mu')}{\d(\mu\otimes\mu')}(x,x')=\frac{\d g_*\mu}{\d\mu}(x)\frac{\d g'_*\mu'}{\d\mu'}(x')\ge CC',$$
    and we are done.
\end{proof}

\subsection{Finiteness of moments}
We fix here a seminorm $\ell$ on $G$, i.e. a (measurable) map $\ell:G\to\m R_{\ge 0}$ such that $\ell(gh)\le \ell(g)+\ell(h)$ for all $g,h\in G$. Assume moreover that $\ell$ is bounded on compact sets. One way to construct such objects is to let $G$ act continuously on some metric space $X$ by isometries, choosing some base point $o\in X$ and letting $\ell(g)=d(o,go)$. Also, if the group $G$ is compactly generated and $K\subset G$ is a compact generating set, one can let $\ell(g)$ be the minimal $n$ such that $g\in K^n$. 

Such a choice of $\ell$ will be implicit in what follows. Note that by restriction, this choice also defines a choice of seminorm on $\Gamma$, which we keep calling $\ell$.

\begin{dfn}
    Let $\mu\in\Prob(G)$. If $k\in\m N$, we say $\mu$ has \emph{finite}  $k$\emph{-th moment} if 
    $$\m E_1(\ell(Z_1)^k)=\int_G \ell(g)^k\d\mu(g)<\infty.$$
    We say it has a \emph{finite exponential moment} if for some $a>0$,
    $$\m E_1(e^{a\ell(Z_1)})=\int_G \exp(a\ell(g))\d\mu(g)<\infty.$$
\end{dfn}

We start with an estimate of moments for convolutions of measures, ensuring that the finiteness of moments of $\mu$ passes down to $\hat\mu$.

\begin{lem}\label{lem: estimations de moments pour Zn et mu^}
    Let $m_k=\m E_1(\ell(Z_1)^k)$ be the $k$-th moment of $\mu$. Then $\m E_1(\ell(Z_n)^k)\le n^km_k$. 
    
    In particular, the $k$-th moment of $\hat\mu$ is bounded by $C_km_k$, where $C_k=\sum_n n^k2^{-n}<\infty$.

    Similarly, if $M_a=\m E_1(e^{a\ell(Z_1)})$ is an exponential moment of $\mu$, then $\m E_1(e^{a\ell(Z_n)})\le M_a^n$. In particular, if $\mu$ has an exponential moment, so does $\hat\mu$.
\end{lem}

\begin{proof}
    We have 
    \begin{eqnarray*}
        \m E_1(\ell(Z_n)^k)&=&\int_{G^n}\ell(g_1\cdots g_n)^k\d\mu^{\otimes n}(g_1,\ldots ,g_n)\\
        &\le&\int_{G^n}\Big(\ell(g_1)+\cdots+\ell(g_n)\Big)^k\d\mu^{\otimes n}(g_1,\ldots ,g_n)\\
        &\le&\int_{G^n}n^{k-1}\Big(\ell(g_1)^k+\cdots+\ell(g_n)^k\Big)\d\mu^{\otimes n}(g_1,\ldots,g_n)\\
        &=&n^{k-1}\sum_{i=1}^n\int_{G^n}\ell(g_i)^k\d\mu^{\otimes n}(g_1,\ldots,g_n)\\
        &=&n^{k-1}\sum_{i=1}^n\int_G\ell(g_i)^k\d\mu(g_i)=n^km_k,
    \end{eqnarray*}
    where the inequality 
    $$\left(\frac{x_1+\cdots +x_n}n\right)^k\le \frac{x_1^k+\cdots+x_n^k}n$$
    is Jensen's inequality for the convex function $x^k$.

    For the exponential part,
    \begin{eqnarray*}
        \m E_1(e^{a\ell(Z_n)})&=&\int_{G^n}\exp(a\ell(g_1\cdots g_n))\d\mu^{\otimes n}(g_1,\ldots,g_n)\\
        &\le&\int_{G^n}\exp(a(\ell(g_1)+\cdots+\ell(g_n)))\d\mu^{\otimes n}(g_1,\ldots,g_n)\\
        &=&\prod_{i=1}^n\int_G e^{a\ell(g_i)}\d\mu(g_i)=M_a^n.
    \end{eqnarray*}
    The $k$-th moment estimate of $\hat\mu$ is now immediate, and the finiteness of an exponential moment will be true as soon as there is some $a>0$ for which $M_a<2$. This comes from the fact that $M_0=1$, on some neighbourhood of $0$ in $\m R_{\ge 0}$, $M_a$ is finite and the map $a\mapsto M_a$ is continuous on that neighbourhood (the integrand is dominated by the value at $a$, which is integrable by assumption).
\end{proof}

We are interested in proving stability results of the finiteness of moments when passing from the measure $\mu$ to the measure $\nu$ in Theorem \ref{theo: BPF lattice in lcsc group}. For the remainder of this section, we fix thus a couple objects that appear in the proof of Theorem \ref{theo: BPF lattice in lcsc group}: the group $G$ and the lattice $\Gamma<G$, the measures $\mu$ and $\nu$, the relatively compact open set $U$ and the constant $C\in]0,1[$, with $U$ symmetric and such that $U\cap\gamma U=\emptyset$ for all $\gamma\in\Gamma$. Finally, we will have to make the additional assumption that $\mu\ll \lambda$.\bigskip

We now define some more stopping times: let $T_n$ be defined as the $n$-th time of hitting $\Gamma U$, or with a formula, $T_{-1}=-1$ and for any $n\ge 0$,
$$T_n=\inf\{k>T_{n-1},Z_k\in\Gamma U\}.$$
For a general measurable set $A\subset \Gamma U$, let also $\tau_A$ be the first return time to $A$, i.e. $\tau_A=\inf\{k\ge 1,Z_k\in A\}$, and denote $\tau=\tau_{\Gamma U}$, i.e. $\tau=T_0\m 1_{\{Z_0\not\in\Gamma U\}}+T_1\m 1_{\{Z_0\in\Gamma U\}}$. Let also $\sigma_A$ be the first hitting time to $A$, that is $\sigma_A=\inf\{k\ge 0,Z_k\in A\}$. \medskip

Recall additionally the stopping times $R_n,S_n,N,\tilde N$ used in Section \ref{ssec: stopping times everywhere}.\medskip

The stopping times we just introduced will be more practical to establish moment estimates with, whereas the stopping times we had introduced before were adapted to our discretization procedure. The following allows us to relate the two.

\begin{lem}
    For any $n\ge 1$, we have $R_n\le T_{2n}$.
\end{lem}

\begin{proof}
    We proceed by induction; note that $S_0\le 1\le T_1$. $R_1$ is the first time the walk hits $\Gamma U$ after $S_0$, hence $R_1\le T_2$. Then again $S_1=R_1+1\le T_2+1\le T_3$ and so on: if for some $k\ge 1$ it holds that $R_k\le T_{2k}$, then $S_k\le 1+T_{2k}\le T_{2k+1}$, and then $R_{k+1}\le T_{2k+2}$ as $R_{k+1}$ is the first time the walk hits $\Gamma U$ after $S_k$.
\end{proof}

Recall that according to Lemma \ref{lem: distri de ZRN est nu}, in order to understand moments of $\nu$, we need only understand the moments of the distribution of the stopped random walk $Z_{R_N}$. 

We start with the following:

\begin{lem}\label{lem: exp moments of Tn}
    If $\Gamma$ is uniform, there are constants $c>1$ and $\alpha>1$ such that 
    $$\m E_1(c^{T_n})\le\alpha^n.$$ 
    Moreover, up to changing $c$, the constant $\alpha$ can be taken arbitrarily close to $1$.
\end{lem}

\begin{proof}
    Note that since $\Gamma$ is uniform, $G/\Gamma$ is compact. 
    
    It is a well known fact that $x\mapsto \mu(xA)$ is continuous for any Borel $A$ and probability measure $\mu\ll \lambda$ on $G$. From there, we may apply successively Theorems 6.0.1, 16.2.5 and finally 16.0.2 ($vii$) of \cite{MT09}, to get that there exists a $c>1$ for which
    \begin{eqnarray}
        \alpha:=\sup_{x\in G/\Gamma}\m E_x(c^\tau)<\infty.\label{eqn: controle uniforme du moment exp du tps de retour}
    \end{eqnarray}

    Note that since $1\in U$, $\m E_1(c^{T_0})=1$. Moreover,
    \begin{eqnarray*}
        \m E_1(c^{T_{n+1}})&=&\m E_1(c^{T_n+\tau\circ\theta_{T_n}})\\
        &=&\m E_1(\m E_1(c^{T_n} c^{\tau\circ\theta_{T_n}}|T_n))\\
        &=&\m E_1(c^{T_n}\m E_1(c^{\tau\circ\theta_{T_n}}|T_n))\\
        &=&\m E_1(c^{T_n}\m E_{Z_{T_n}}(c^{\tau}))\\
        &\le&\alpha\m E_1(c^{T_n}) \qquad\text{by (\ref{eqn: controle uniforme du moment exp du tps de retour}),}
    \end{eqnarray*}
    and the first claim follows by induction.\medskip

    For a fixed $x$, note that the map
    $$t\mapsto \m E_x(t^\tau)=\sum_{k\ge 0}\m P_x(\tau=k)t^k$$
    is convex on $\m R_{\ge 0}$. As such, the map 
    $$t\mapsto \sup_{x\in G/\Gamma}\m E_x(t^\tau)$$
    is also convex, and it is in particular continuous on $\m R_{>0}$. The second claim follows.
\end{proof}

\begin{lem}\label{lem: poly moments of Tn}
    If $\Gamma$ is uniform, for any $k\in\m N$ there is a constant $\beta>0$ such that 
    $$\m E_1({T_n}^k)\le\beta n^k.$$
\end{lem}

\begin{proof}
    The proof is essentially the same as that of Lemma \ref{lem: exp moments of Tn}, the biggest difference being the use of the Minkowski inequality.\medskip
    
    Let
    \begin{eqnarray}
        \beta:=\sup_{x\in G/\Gamma}\m E_x(\tau^k),
    \end{eqnarray}
    which is finite by (\ref{eqn: controle uniforme du moment exp du tps de retour})

    Note that since $1\in U$, $\m E_1({T_0}^k)=0$. Moreover,
    \begin{eqnarray*}
        \m E_1({T_{n+1}}^k)^{1/k}&=&\m E_1((T_n+\tau\circ\theta_{T_n})^k)^{1/k}\\
        &\le&\m E_1(T_n^k)^{1/k} +\m E_1(\m E_1(\tau^k\circ\theta_{T_n}|T_n))^{1/k}\\
        &=&\m E_1({T_n}^k)^{1/k}+\m E_1(\m E_{Z_{T_n}}(\tau^k))^{1/k}\\
        &\le&\m E_1({T_n}^k)^{1/k}+\beta^{1/k}
    \end{eqnarray*}
    and the claim follows by induction.
\end{proof}

We now set out to prove that even when $\Gamma$ is non-uniform, the return time is uniformly integrable on compact sets. In a discrete setting this would fall out of the positive recurrence of the chain, but the non-discrete setting adds some technical difficulties. We start with the following technical lemma.

\begin{lem}\label{lem: continuity of proba of return time}
    The map $x\mapsto \m E_x(\tau_A)$ is continuous.
\end{lem}

\begin{proof}
    Let $f:x\mapsto \m E_x(\sigma_A)$, and notice that $\m E_x(\tau_A)=Pf(x)$. By \cite[Theorem 6.1.1]{MT09}, we need only show that for any measurable set $B$, the map $P\m 1_B$ is continuous.

    Since $P\m 1_B(x)=\mu(x\inv B)$, continuity follows from the fact that $\mu\ll\lambda$ and that $\mu$ has finite measure.
\end{proof}

\begin{lem}
    For any relatively compact set $K\subset G/\Gamma$ with $m(K)>0$, we have
    $$\sup_{x\in K}\m E_x(\tau_K)<\infty$$
\end{lem}

\begin{proof}
    Recall that a \emph{regular set} is a set $A$ such that for any $B$ with $m(B)>0$,
    $$\sup_{x\in A}\m E_x(\tau_B)<\infty.$$
    Our aim is to show that $K$ is regular, for which we first show that $(\pi(Z_n))_n$ is regular, that is, $G/\Gamma$ can be covered by countably many regular sets.

    By Proposition \ref{pro: lattice in group is Harris}, $(\pi(Z_n))_n$ is positive Harris recurrent, so by \cite[Theorems 11.0.1, 11.3.11]{MT09}, there is a set $A$, some $b>0$ and a map $V:G/\Gamma\to\m R_{\ge 0}$ such that 
    $$PV-V\le -1+b\m 1_A,$$
    and $x\mapsto\m E_x(\tau_A)$ is bounded on $A$ and finite on $\{V<\infty\}$. Say 
    $$\sup_{x\in A}\m E_x(\tau_A)=M.$$
    We may apply \cite[Theorem 16.2.6 $(ii)\implies (d)$]{DMPS18} with $f=r\equiv 1$ and $V_n=V$ to get that $\{V<\infty\}$ is the set of regular points in $G/\Gamma$, i.e. the set of points $x$ for which for any set $B$ with $m(B)>0$, $\m E_x(\tau_B)<\infty$. However, the set of regular points is both non-empty and $G$-invariant, and as such $\{V<\infty\}=G/\Gamma$, and $\m E_x(\tau_A)$ is finite everywhere.
    
    By \cite[Proposition 11.3.13]{MT09}, the level sets 
    $$C_\ell=\{x\in G/\Gamma,\m E_x(\sigma_A)\le \ell\}$$
    are all regular, and since $\sigma_A=\tau_A\m 1_{\{Z_0\in A\}}$, we have $\m E_x(\sigma_A)=\m E_x(\tau_A)$ for all $x\not\in A$. Hence, for all $\ell\ge M$, we also have 
    $$C_\ell=\{x\in G/\Gamma,\m E_x(\tau_A)\le \ell\}.$$
    Moreover, by Lemma \ref{lem: continuity of proba of return time}, $\m E_x(\tau_A)$ is continuous, and in particular the sublevel sets $C_\ell$ are closed for $\ell\ge M$. Moreover, 
    $$G/\Gamma=\bigcup_{\ell\ge M}C_\ell$$
    and $G/\Gamma$ is Polish, so by Baire's category theorem, there is some $\ell$ for which $C_\ell$ has non-empty interior. Since subsets of regular sets clearly are regular, we have some open set $W$ which is regular. Note that
    $$\overline K\subset G/\Gamma=\bigcup_{g\in G}gW,$$
    so by compactness of $\overline K$ there are some $g_1,\dots, g_n$ for which 
    $$K\subset \bigcup_{i=1}^ng_iW.$$
    Regularity is preserved under the $G$-action and finite unions, so $K$ is regular and we are done.
\end{proof}

\begin{lem}\label{lem: Tn est sslinéaire}
    There is a constant $\beta$ such that for all $n\in\m N$, 
    $$\m E_1(T_n)\le\beta n.$$
\end{lem}

\begin{proof}
    Let
    $$\beta=\sup_{x\in \Gamma U}\m E_x(\tau),$$
    which is finite by the previous lemma. Now 
    \begin{eqnarray*}
        \m E_1(T_{n+1})&=&\m E_1(T_n+\tau\circ\theta_{T_n})\\
        &=&\m E_1(T_n)+\m E_1(\m E_1(\tau\circ\theta_{T_n}|T_n))\\
        &=&\m E_1(T_n)+\m E_1(\m E_{Z_{T_n}}(\tau))\\
        &\le&\m E_1(T_n)+\beta
    \end{eqnarray*}
    since $Z_{T_n}\in\Gamma U$ $\m P_1$-a.s. Since $1\in U,$ we have $\m E_1(T_0)=0$ and the result follows by induction.
\end{proof}

Now that we have integrability results of the return times, we can use them to get integrability results for $\nu$.

\begin{pro}\label{pro: moments polynomiaux de nu}
    Let $k\in\m N$. Assume that $\mu$ has finite $k$-th moment and $\m E_1((T_n)^k)\le \beta n^k$ for some $\beta>0$. Then $\nu$ has finite $k$-th moment.
\end{pro}

\begin{proof}
    Let $m_k=\m E_1(\ell(Z_1)^k)<\infty$. By Lemma \ref{lem: estimations de moments pour Zn et mu^}, we have $\m E_1(\ell(Z_n)^k)\le~m_kn^k$ and 
    \begin{eqnarray*}
        \tilde{\m E}_1(\ell(Z_{R_N})^k)&\le&\tilde{\m E}_1(\ell(Z_{T_{2N}})^k)\\
        &=&\tilde{\m E}_1(\tilde{\m E}_1(\ell(Z_{T_{2N}})^k|T_{2N}))\\
        &\le&\tilde{\m E}_1(m_k(T_{2N})^k)\\
        &=&m_k\tilde{\m E}_1(\tilde{\m E}_1((T_{2N})^k|N))\\
        &\le&m_k\tilde{\m E}_1(\beta (2N)^k)\le m_k\tilde{\m E}_1(\beta (2\tilde N)^k)<\infty
    \end{eqnarray*}
    since geometric laws have finite polynomial moments.
    
    Let now $A$ be the (finite) diameter of the relatively compact set $U$, so that for all $\gamma\in\Gamma$ and $x\in\gamma U,$
    $$\ell(\gamma)\le\ell(x)+A.$$
    But then,
    $$\sum_{\gamma\in\Gamma}\ell(\gamma)^k\nu(\gamma)=\sum_{\gamma\in\Gamma}\ell(\gamma)^k\tilde{\m P}_1(Z_{R_N}\in\gamma U)\le \tilde{\m E}_1((\ell(Z_{R_N})+A)^k)<\infty$$
    as desired.
\end{proof}

\begin{pro}\label{pro: moments exponentiels de nu}
    Assume that $\mu$ has a finite exponential moment, and that for some $\alpha\in]1,(1-C^2)^{-1/2}[$, there is some $c>1$ such that $\m E_1(c^{T_n})\le \alpha^n$. Then $\nu$ has a finite exponential moment.
\end{pro}

\begin{proof}
    Since $\mu$ has a finite exponential moment, the same proof as Lemma \ref{lem: estimations de moments pour Zn et mu^} shows that for some $a>0$, $M_a=\m E_1(e^{a\ell(Z_1)})\le c$, so $\m E_1(e^{a\ell(Z_n)})\le c^n$. Now
    \begin{eqnarray*}
        \tilde{\m E}_1\left(e^{a\ell(Z_{R_N})}\right)&\le&\tilde{\m E}_1\left(e^{a\ell(Z_{T_{2N}})}\right)\\
        &=&\tilde{\m E}_1\left(\tilde{\m E}_1\left(e^{a\ell(Z_{T_{2N}})}|T_{2N}\right)\right)\\
        &\le&\tilde{\m E}_1\left(c^{T_{2N}}\right)\\
        &=&\tilde{\m E}_1\left(\tilde{\m E}_1\left(c^{T_{2N}}|N\right)\right)\\
        &\le&\tilde{\m E}_1\left(\alpha^{2N}\right)\le \tilde{\m E}_1\left(\alpha^{2\tilde N}\right)<\infty
    \end{eqnarray*}
    since $\alpha^2<\frac1{1-C^2}$. Finally, if $A$ denotes again the diameter of $U$, we have
    $$\sum_{\gamma\in\Gamma}e^{a\ell(\gamma)}\nu(\gamma)= \sum_{\gamma\in\Gamma}e^{a\ell(\gamma)}\tilde{\m P}_1(Z_{R_N}\in\gamma U)\le\tilde{\m E}_1\left(e^{a(\ell(Z_{R_N})+A)}\right)<\infty$$
    and we are done.
\end{proof}

We finish this section by compiling the various results we have just proven in the following:

\begin{theo}\label{theo: moment preservation}
    Let $\mu$ be as in Theorem \ref{theo: BPF lattice in lcsc group} and assume moreover that $\mu\ll\lambda$, where $\lambda$ is the Haar measure on $G$. If $\mu$ has finite first moment, then so does $\nu$.\medskip

    If now $\Gamma$ is uniform and $\mu$ has finite $k$-th moment, then so does $\nu$. If $\mu$ has a finite exponential moment, then so does $\nu$.
\end{theo}

\begin{proof}
    From Lemma \ref{lem: Tn est sslinéaire}, we have a constant $\beta$ such that $\m E_1(T_n)\le \beta n$, so we may apply Proposition \ref{pro: moments polynomiaux de nu} with $k=1$ to get the first part of the theorem.\medskip

    Assume now that $\Gamma$ is uniform and that $\mu$ has finite $k$-th moment. Then by Lemma \ref{lem: poly moments of Tn}, there is a constant $\beta$ such that $\m E_1({T_n}^k)\le \beta n^k$. We may thus apply Proposition \ref{pro: moments polynomiaux de nu} to conclude.

    If now $\mu$ has a finite exponential moment, we may apply Lemma \ref{lem: exp moments of Tn} to get a constant $c$ such that $\m E_1(c^{T_n})\le \alpha^n$ for some $\alpha$ with $1<\alpha<\frac1{\sqrt{1-C^2}}$. Proposition \ref{pro: moments exponentiels de nu} now ensures that $\nu$ also has a finite exponential moment.
\end{proof}

\subsection{The case where $G$ is a product}
We end this section with a proof of Corollary \ref{corintro: BPF de réseau dans produit}.

\begin{proof}[Proof of Corollary \ref{corintro: BPF de réseau dans produit}]
    By Lemma \ref{lem: stabilité de (H) sous produit} and the fact that all $\mu_i$ have (H), the measure $\mu=\mu_1\otimes \mu_2\otimes\cdots\otimes \mu_n$ on $G$ has (H) as well. Now Theorem \ref{theointro: LS pour des marches sur des groupes lcsc} gives a probability measure $\nu$ on $\Gamma$ such that 
    $$B_{PF}(\Gamma,\nu)\simeq B_{PF}(G,\mu).$$
    That 
    $$B_{PF}(G,\mu)\simeq B_{PF}(G_1,\mu_1)\times\cdots\times B_{PF}(G_n,\mu_n)$$
    is a consequence of \cite[Corollary 3.2]{BS06}. Combining the two isomorphisms, we get the desired result.
\end{proof}

\section{Pulling back walks from graphs to groups}\label{sec: kaima woess+recurrence}

In \cite{Der25}, we gave a way to transform a walk on a graph to a walk on a \emph{lattice} on that graph, that is, a discrete group acting on the graph with ``finite covolume", and this construction preserved the Poisson-Furstenberg boundary. We wish to slightly generalize this result here, removing the discreteness hypothesis on the group (but keeping the cofiniteness hypothesis). We wish to prove this generalization in view of the following reasoning: \medskip

Assume that $\Gamma<G\times H$ is a lattice and that $H$ is acting on a graph $X$ with finite covolume. Then using Corollary \ref{corintro: BPF de réseau dans produit}, for decent choices of probability measures $\mu_G$ and $\mu_H$ on $G$ and $H$ we have a measure $\nu$ on $\Gamma$ such that 
$$B_{PF}(\Gamma,\nu)\simeq B_{PF}(G,\mu_G)\times B_{PF}(H,\mu_H)$$
If moreover it happens that the Poisson boundary of the graph $X$ for some random walk $P$ on it is well understood, one could hope that there is a measure $\mu_H$ on $H$ which remembers $P$, so that 
$$B_{PF}(\Gamma,\nu)\simeq B_{PF}(G,\mu_G)\times B_{PF}(H,\mu_H)\simeq B_{PF}(G,\mu_G)\times B_{PF}(X,P).$$
Hence if the Poisson boundaries of $\mu_G$ and $P$ are well understood, a Poisson boundary of $\Gamma$ will be well understood as well, thanks to Theorem \ref{theo: RW on graph -> on group}.

\begin{dfn}\label{dfn: covolume/réseau dans un graphe}
    Let $X$ be a connected, locally finite graph, we endow $\Aut(X)$ with the compact open topology. A closed \emph{unimodular} subgroup $G<\Aut(X)$ is said to be \emph{of finite covolume} in $X$ if for some fundamental domain $D$ for the action of $G$ on (the vertices of) $X$,
    $$\sum_{x\in D}\frac1{\lambda(G_x)}<\infty,$$
    where $\lambda$ is a left Haar measure of $G$. This holds in particular if $G$ is cocompact. Note that since the stabilizer $G_x$ is compact and open, $0<\lambda(G_x)<\infty$ and the value of the sum does not depend on the choice of a fundamental domain by unimodularity.
\end{dfn}

If moreover $G$ is discrete in $\Aut(X)$ (equivalently, if every point stabilizer is finite), $G$ is said to be a \emph{lattice} in $X$. 

\begin{theo}\label{theo: RW on graph -> on group}
    Let $X$ be a locally finite connected graph of bounded degree, and let $G<\Aut(X)$ be closed. Let $P$ be the Markov operator associated to the simple random walk. If $G$ is either acting cocompactly or is unimodular of finite covolume, there exists a measure $\mu$ on $G$ for which 
    $$B_{PF}(X,P)\simeq B_{PF}(G,\mu)$$
    $G$-equivariantly, and $\mu$ can be chosen to be symmetric if $X$ has constant degree and $G$ is unimodular.
\end{theo}

The strategy of proof is the same as that of \cite[Theorem 1]{Der25}. More specifically, we first prove that an orbit $G\cdot o$ of $G$ in $X$ is recurrent for the random walk, in the sense that the walk meets the orbit infinitely often almost surely. Equivalently, we prove that the quotient walk on $G\backslash X$ is recurrent. We may then use \cite[Theorem 3.2]{Der25} to get a $G$-invariant random walk on $G\cdot o$ which has the same bounded harmonic functions as those on $X$. The last step is to apply Propositions 2.15 and 3.1 of \cite{KW02} to pull the walk on $G\cdot o$ back to a walk on $G$.

\begin{pro}\label{pro: (positive) recurrence criterion}
    If $X$ is of bounded degree and if $G$ is unimodular and of finite covolume in $X$ and $o$ is any vertex, the simple random walk on $X$ meets the orbit $G\cdot o$ infinitely often, regardless of the starting distribution.
\end{pro}

The proof of this result is the same as that of \cite[Corollary 3.6]{Der25}. In fact, the network version of this result remains valid in the context of unimodular locally compact groups.

\begin{proof}
    Consider the quotient $G\backslash X$ endowed with the quotient walk 
    $$\bar p(\bar x,\bar y)=\sum_{y\in\bar y}p(x,y),$$ 
    where the sum does not depend on the choice of the lift $x$ of $\bar x$. We prove that the measure $m(\bar x)=\frac{\deg(x)}{\lambda(G_x)}$ (which does not depend on the choice of the lift) is stationary. This will ensure that the walk on the quotient is (positive) recurrent by Proposition \ref{pro: inv proba->recur}, as $m$ is finite from the assumptions on $X$ and $G$. 

    For stationarity, let $x\in X$. Then
    \begin{eqnarray*}
        m\bar P(\bar x)
        &=&\sum_{\bar y\in G\backslash X}\frac{\deg \bar y}{\lambda(G_{\bar y})}\sum_{x\in\bar x}p(y,x)\\
        &=&\sum_{\bar y\in G\backslash X}\frac{\deg y}{\lambda(G_y)}\frac1{\lambda(G_x)}\int_G p(y,gx)\d\lambda(g)\\
        &\overset{(1)}=&\sum_{\bar y\in G\backslash X}\frac{\deg y}{\lambda(G_x)}\frac1{\lambda(G_y)}\int_G p(g\inv y,x)\d\lambda(g)\\
        &\overset{(2)}=&\sum_{\bar y\in G\backslash X}\frac{\deg y}{\lambda(G_x)}\frac1{\lambda(G_y)}\int_G p(gy,x)\d\lambda(g)\\
        &=&\sum_{\bar y\in G\backslash X}\frac{\deg y}{\lambda(G_x)}\sum_{y\in \bar y}p(y,x)\\
        &=&\sum_{y\in X}\frac1{\lambda(G_x)}\m 1_{x\sim y}(x,y)=\frac{\deg(x)}{\lambda(G_x)}=m(\bar x),
    \end{eqnarray*}
    where $(1)$ is the fact that $p$ is $G$-invariant, and $(2)$ comes from unimodularity of $G$. This concludes the proof.
\end{proof}

\begin{rem}
    In much the same way as \cite[Proposition 3.7]{Der25}, for a given network $\mc N$ it is possible to define a quotient network $G\backslash \mc N$ for any unimodular closed subgroup $G<\Aut(\mc N)$ by the same formula 
    $$\bar a(\bar x,\bar y)=\frac1{\lambda(G_x)}\sum_{y\in \bar y}a(x,y)$$
    for some lift $x$ of $\bar x$. The total conductance $m'(x)=m(x)/\lambda(G_x)$ is again a stationary measure on $G\backslash X$, and its finiteness implies positive recurrence of the quotient walk, in the same way as Proposition \ref{pro: (positive) recurrence criterion} and \cite[Proposition 3.5]{Der25}
\end{rem}

\begin{proof}[Proof of Theorem \ref{theo: RW on graph -> on group}]
    Let $o\in X$ be fixed. If $G$ is unimodular and of finite covolume, Proposition \ref{pro: (positive) recurrence criterion} tells us that the orbit $G\cdot o$ is recurrent. If $G$ is acting cocompactly, the quotient walk is irreducible on a finite state space, hence recurrent. In either case, we may apply \cite[Theorem 3.2]{Der25} to say that there is a $G$-invariant random walk $Q$ on $G\cdot o$ for which 
    $$B_{PF}(X,P)=B_{PF}(G\cdot o,Q)$$
    $G$-equivariantly, and $Q$ is symmetric as soon as $P$ is. We may now apply the correspondence of \cite[Proposition 2.15]{KW02} to say that $Q$ corresponds to the walk $\mu$ on $G$ defined by 
    $$\frac{\d\mu}{\d\lambda}(g)=q(o,go),$$
    where $\lambda$ is the Haar measure on $G$, normalized so that $\lambda(G_o)=1$. We may finally apply \cite[Proposition 3.1]{KW02} to conclude that the Poisson boundaries of $(G\cdot o,Q)$ and $(G,\mu)$ are ($G$-equivariantly) isomorphic. It is easily checked that if $G$ is unimodular and $Q$ is symmetric, then so is $\mu$.
\end{proof}

One noteworthy remark in the scheme presented at the beginning of this section is that we may apply results from Section \ref{sec: discretisation} to the measure we just constructed.

\begin{pro}\label{pro: pullback de RW de graphe a (H)}
    The measure $\mu$ from Theorem \ref{theo: RW on graph -> on group} satisfies condition (H).
\end{pro}

\begin{proof}
    Let $k$ be in the stabiliser $G_o$ of $o$. Note that 
    $$p(o,kgo)=p(k\inv o,go)=p(o,go),$$ 
    so $k_*\mu=\mu$.
    However, $G_o$ is open, so we may take $U=G_o$ and $C$ arbitrarily close to $1$ in the definition of  condition (H) (see Definition \ref{dfn: condition (H)}).
\end{proof}

\begin{rem}
    The same results on finiteness of moments stated in \cite[Proposition 3.12]{Der25} hold for this walk as well, with the exact same proofs.
\end{rem}

\bibliographystyle{alpha}  
\bibliography{biblio}
\end{document}